\documentclass[a4paper]{article}
\usepackage{amsmath,amssymb,fancybox,longtable,graphics}
\usepackage[dvipdfm]{graphicx}
\newtheorem{thm}{\textbf{Theorem}}[section]
\newtheorem{proposition}[thm]{\textbf{Proposition}}

\def\QED{$\Box$} 
\def\mbi#1{\boldsymbol{#1}} 
\def\Spec{\mathop{\mathrm{Spec}}\nolimits}
\def\Proj{\mathop{\mathrm{Proj}}\nolimits}
\def\Pic{\mathop{\mathrm{Pic}}\nolimits}

\begin{document}
\title{Isomorphism among families of weighted $K3$ hypersurfaces}
\author{
\textsc{Masanori} KOBAYASHI \\
{\small{\textit{Department of Mathematics and Information Sciences, Tokyo Metropolitan University, }}}\\
{\small{\textit{1-1 Minami-Osawa, Hachioji-shi Tokyo, 192-0397, Japan }}}\\
{\small{\textit{kobayashi-masanori@tmu.ac.jp }}}\\
\\
\textsc{Makiko} MASE \\
{\small{\textit{Department of Mathematics and Information Sciences, Tokyo Metropolitan University, }}}\\
{\small{\textit{1-1 Minami-Osawa, Hachioji-shi Tokyo, 192-0397, Japan }}}\\
{\small{\textit{mase-makiko@ed.tmu.ac.jp } }}}
\date{\null}
\maketitle

Some of the 95 families of weighted $K3$ hypersurfaces have been known to have the isometric lattice polarizations. 
It is shown that weighted $K3$ hypersurfaces in such families are to one-to-one correspond by explicitly constructing the monomial birational morphisms among the weighted projective spaces. 
All the weight systems having the isometric Picard lattices commonly possess an anticanonical sublinear system, being confirmed that the Picard lattice of the sublinear system we obtained is the same as those of the complete linear systems. 

\vspace{3mm}
\noindent
Mathematics Subject Classification 2000:  14J28, 14J10, 14J17.

\section{Introduction}

There is a famous list of 95 families of weighted $K3$ hypersurfaces \cite{Yonemura}\cite{Fletcher}. 
The Picard lattice of the minimal model of a generic member of each family has been computed by Belcastro \cite{Belcastro}. 
Here the Picard lattice of a $K3$ surface $S$ means the 
Picard group $\Pic S$ together with the cup product. 
Some of the lattices are found to be isometric. 
So one may expect that the corresponding families of $K3$ surfaces would coincide, 
in the sense that the period maps have the same image. 
It is impossible to identify the whole complete anticanonical linear systems, 
since the dimensions of the systems do not always coincide. 
Nevertheless, we show that there exists an identification between subfamilies. In fact, for all such pairs, there exists an explicit monomial birational map between the weighted projective spaces which induces, on the families of minimal models, an isomorphim between subfamilies of $K3$ surfaces. 
Moreover, the subfamilies are general enough, namely, the generic member of each subfamily has the same Picard lattice as the original family. 
We remark that these maps are compatible with logarithmic moment maps and keep the amoebas of $K3$ surfaces. 

%
%
%

We first explain the idea by a toy model, 
namely the family of elliptic curves in the projective plane. 
Let $P_1$ and $P_2$ be distinct points in $\mbi{P}^2$, and $L$ be the line through them. 
Blow up $P_1$ and $P_2$, blow down the strict transform of $L$, 
and we get $\mbi{P}^1 \times \mbi{P}^1$. 
We write the transform of $P_1,P_2$ and $L$ in $\mbi{P}^1 \times \mbi{P}^1$ as $H_1,H_2$ and $Q$, respectively. 
A general cubic in $\mbi{P}^2$ is mapped to an element of $|3H_1+3H_2-3Q|$ in $\mbi{P}^1 \times \mbi{P}^1$, 
which is not anticanonical. 
Take an anticanonical sublinear system $\mbi{L}=|3L-P_1-P_2|$, 
which is seven-dimensional. 
The transform of $\mbi{L}$ in $\mbi{P}^1 \times \mbi{P}^1$ is $|2H_1 + 2H_2 - Q|$, 
which is an anticanonical sublinear system. 
The correspondence of Newton polygons is discribed below, 
which also involves the complete anticanonical linear system of the Del Pezzo surface of degree $7$. 

\begin{figure}[htbp]
\begin{center}
\includegraphics[width=.30\linewidth]{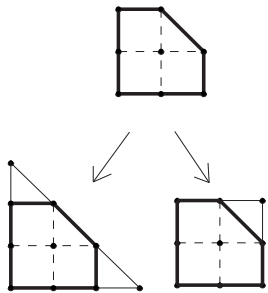}
\end{center}
\caption{Correspondence of Newton polygons.\label{fig. 1. }}
\end{figure}

\section{Setup}

Let $a\, := (a_0,a_1,a_2,a_3)$ be a list of positive integers, which are called {\it weights.\/} 
Let $\mbi{P}(a)$ be the weighted projective space $\Proj \mbi{C}[W,X,Y,Z]$ 
where degrees of $W,X,Y,Z$ are $a_0,a_1,a_2, a_3$, respectively. 
We can assume without a loss of generality that 
$a_0 \leq a_1 \leq a_2 \leq a_3$ and also that 
the weights are {\it well-posed}, that is, every greatest common divisor of all but one of the $a_i$'s is one. 
Let $M(a)$ be the group of exponents of degree-zero rational monomials 
\[
\left\{ (m_0, m_1, m_2, m_3) \, \in\, \mbi{Z}^4\, \left| \, \sum_{i=0}^3\, a_i m_i= 0 \right.\right\}.
\] 
It is easy to see that $M := \mbi{Z}^3 \cong M(a)$. 
Define a rational tetrahedron 
\[
\Delta(a):=\{(m_0,m_1,m_2,m_3) \in M(a) \otimes \mbi{R} \ | \ m_i \geq -1\}.
\] 
After the multiplication of the monomial $WXYZ$, 
$\Delta(a) \cap M(a)$ generates sections of the anticanonical bundle of $\mbi{P}(a)$ as a vector space. 

In general, given a bounded rational convex polyhedron $\Delta$ in $\mbi{R}^n$, 
one has an $n$-dimensional 
projective 
toric variety 
$\mbi{P}_{\Delta}$ in a standard way. 
$\mbi{P}_{\Delta(a)}$ is isomorphic to $\mbi{P}(a)$, 
and contains the three-dimensional algebraic torus $\mbi{T}:= \Spec\mbi{C}[M]$. 

A 
convex subpolyhedron $\Delta$ in $\Delta(a)$ determines an anticanonical linear subsystem, which corresponds to a family of Laurent polynomials in $\mbi{C}[M]$ with degree zero and Newton polytope in $\Delta$. 
Each polynomial $F$ determines the zero set $Z_F$ in $\mbi{T}$. 
The zero set
$\overline{Z_F}$ in $\mbi{P}_{\Delta}$ usually has singularities. 
The minimal resolution $S_F$ of $\overline{Z_F}$ is a $K3$ surface if and only if 
$\overline{Z_F}$ has no boundary component and 
the singularities are only rational double points. 
In that case, for a generic $F$ we denote the Picard lattice $\Pic S_F$ by $\Lambda_\Delta$. 
If $\Delta=\Delta(a)$, 
$S_F$ is $K3$ 
if and only if $a$ is in the `famous 95' list. 
In this case we write $\Lambda_\Delta$ as $\Lambda(a)$. 

Assume that $\Delta$ is a bounded lattice polyhedron and contains the origin as the only lattice point in its interior. 
$\Delta$ is called {\it reflexive\/} if the polar dual of $\Delta$ is also a lattice polyhedron\cite{Batyrev}. 
For $3$-dimensional $\Delta$, 
the reflexivity 
is equivalent to that 
$\mbi{P}_\Delta$ is a Fano 3-fold with only canonical Gorenstein singularities, 
and that the minimal model of a general anticanonical member is a $K3$ surface.

\section{Result}
We state the main theorem.
\begin{thm} Let two weights $a$ and $b$ be in the `famous $95$' with isometric 
Picard lattices. 
Then there exist subspaces $D_a$ (resp. $D_b$) of the anticanonical complete linear system of $\mbi{P}(a)$ (resp.$\mbi{P}(b)$), 
and an isomorphism $\varphi: D_a \to D_b$ with the following properties. 
(1) If the minimal model of $X \in D_a$ is a K3 surface, then the minimal model of $\varphi(X)$ is an isomorphic K3 surface 
to $X$ as minimal models, 
and vice versa. 
(2) The Picard lattices of the minimal models of generic members of $D_a,D_b$ are isometric to $\Lambda (a) \simeq \Lambda (b)$. 
\end{thm}

The theorem follows from the proposition below: 
\begin{proposition}
Under the assumption of the theorem, 
there exists a group isomorphism $M(a) \cong M(b)$, 
and a common reflexive subpolyhedron $\Delta$ of $\Delta(a)$ and $\Delta(b)$, 
with the following properties. 
(1) the associated birational maps 
$\, \varphi_a : \mbi{P}_\Delta \text{- -$\to$} \mbi{P}(a)$ and 
$\, \varphi_b: \mbi{P}_\Delta \text{- -$\to$} \mbi{P}(b)$ 
map the general anticanonical members of $\mbi{P}_\Delta$ to those of $\mbi{P}(a)$ and $\mbi{P}(b)$, 
(2) The lattices $\Lambda(a)$, $\Lambda(b)$ and $\Lambda_\Delta$ are isometric.
\end{proposition}
%

{\small
\begin{longtable}{rllc}
\hline
No. & Families & The vertices of $\Delta$ & Picard lattice\\
\hline
\hline
13   & $\mbi{P} (1,3,8,12) \supset (24)$ & $ Z^2, W^{24}, W^3X^7, WX^5Y, Y^3, X^4Z $  & $E_6 \perp U$ \\
72   & $\mbi{P} (1,2,5,7) \supset (15)$  & $ WZ^2, W^{15}, WX^7, X^5Y, Y^3, X^4Z $ & (8)\\
\hline
50   & $\mbi{P} (1,4,10,15) \supset (30)$ & $ Z^2, W^{30}, W^2X^7, X^5Y, Y^3 $ & $E_7\perp U$ \\
82   & $\mbi{P} (1,3,7,11) \supset (22)$  & $ Z^2, W^{22}, WX^7, X^5Y, WY^3 $ & (9)\\
\hline
9    & $\mbi{P} (1,4,5,10) \supset (20)$ & $ W^{20}, X^5, Z^2, Y^2Z, WXY^3, W^5Y^3$ & $T_{2,5,5}$ \\
71   & $\mbi{P} (1,3,4,7) \supset (15)$  & $ W^{15}, X^5, WZ^2, Y^2Z, XY^3, W^3Y^3$ & (10)\\
\hline
14   & $\mbi{P} (1,6,14,21) \supset (42)$ &  $ Z^2, Y^3, X^7, W^{42} $& $E_8 \perp U$ \\
28   & $\mbi{P} (1,3,7,10) \supset (21)$  &  $ WZ^2, Y^3, X^7, W^{21} $& (10)\\
45   & $\mbi{P} (1,4,9,14) \supset (28)$  &  $ Z^2, WY^3, X^7, W^{28} $& \\ 
51   & $\mbi{P} (1,5,12,18) \supset (36)$ &  $ Z^2, Y^3, WX^7, W^{36} $& \\
\hline
38   & $\mbi{P} (1,6,8,15) \supset (30)$  & $ Z^2, W^{30}, X^5, XY^3, W^6Y^3 $ & $E_8\perp A_1\perp U$ \\
77   & $\mbi{P} (1,5,7,13) \supset (26)$  & $ Z^2, W^{26}, WX^5, XY^3, W^5Y^3 $ & (11)\\
\hline
20   & $\mbi{P} (1,6,8,9) \supset (24)$  & $ W^6Z^2, W^{24}, X^4, XZ^2, Y^3 $ & $E_8\perp A_2\perp U$ \\
59   & $\mbi{P} (1,5,7,8) \supset (21)$  & $ W^5Z^2, W^{21}, WX^4,  XZ^2, Y^3 $ & (12)\\
\hline
26   & $\mbi{P} (2,4,5,9) \supset (20)$  & $ WZ^2, W^{10}, X^5, Y^4 $   & $D_8\perp D_4\perp U$ \\
34   & $\mbi{P} (2,6,7,15) \supset (30)$ & $ Z^2, W^{15}, X^5, WY^4 $ & (14)\\
\hline
26   & $\mbi{P} (2,4,5,9) \supset (20)$ & $ WZ^2, W^5Y^2, Y^4, X^5, W^8X $ & $D_8\perp D_4\perp U$ \\
34   & $\mbi{P} (2,6,7,15) \supset (30)$ & $ Z^2, W^8Y^2, WY^4, X^5, W^{12}X  $ & (14)\\
76   & $\mbi{P} (2,5,6,13) \supset (26)$ & $ Z^2, W^8X^2, X^4Y, WY^4, W^{13}$ & \\
\hline
27   & $\mbi{P} (2,3,8,11) \supset (24)$  & $ WZ^2, W^{12}, X^8, Y^3 $ & $E_8\perp D_4\perp U$ \\
49   & $\mbi{P} (2,5,14,21) \supset (42)$ & $ Z^2, W^{21}, WX^8, Y^3 $ & (14)\\
\hline
16 & $\mbi{P} (3,6,7,8) \supset (24)$ & $ Z^3, W^3YZ, W^6X, X^4, WY^3$ & $E_8 \perp (A_2)^3 \perp U$ \\
54 & $\mbi{P} (3,5,6,7) \supset (21)$ & $Z^3, W^3XZ, W^7, WY^3, X^3Y$ & (16)\\
\hline
43   & $\mbi{P} (3,4,11,18) \supset (36)$ & $ Z^2, W^{12}, X^9, WY^3 $ & $E_8\perp E_6\perp U$ \\
48   & $\mbi{P} (3,5,16,24) \supset (48)$ & $ Z^2, W^{16}, WX^9, Y^3 $  & (16)\\
\hline
43   & $\mbi{P} (3,4,11,18) \supset (36)$ & $ Z^2, W^6Z, W^8X^3, X^9, WY^3, W^7XY $ & $E_8\perp E_6\perp U$ \\
48   & $\mbi{P} (3,5,16,24) \supset (48)$ & $ Z^2, W^8Z, W^{11}X^3, WX^9, Y^3, W^9XY $  & (16)\\
88   & $\mbi{P} (2,5,9,11) \supset (27)$  & $ XZ^2, W^8Z, W^{11}X, WX^5, Y^3, W^9Y $ & \\
\hline
68   & $\mbi{P} (3,4,10,13) \supset (30)$ & $ XZ^2, X^5Y,  W^2X^6,  Y^3, W^{10} $ & $E_8\perp E_7\perp U$ \\
83   & $\mbi{P} (4,5,18,27) \supset (54)$ & $ Z^2, W^9Y, W^{11}X^2, Y^3, WX^{10} $ & (17)\\
92   & $\mbi{P} (3,5,11,19) \supset (38)$ & $ Z^2, W^9Y, W^{11}X, XY^3, WX^7 $ & \\
\hline
30   & $\mbi{P} (5,7,8,20) \supset (40)$ & $ Z^2, W^4Z, WX^5, W^5XY, Y^5 $   & $E_8\perp T_{2,5,5}$ \\
86   & $\mbi{P} (4,5,7,9) \supset (25)$  & $ YZ^2, W^4Z, X^5, W^5X, WY^3 $  & (18)\\
\hline
46   & $\mbi{P} (5,6,22,33) \supset (66)$ &  $ Z^2, W^{12}X, X^{11}, Y^3 $ & $E_8^2\perp U$ \\
65   & $\mbi{P} (3,5,11,14) \supset (33)$ &  $ XZ^2, W^{11}, WX^6, Y^3 $ & (18)\\
80   & $\mbi{P} (4,5,13,22) \supset (44)$ &  $ Z^2, W^{11}, WX^8, XY^3 $ & \\
\hline
56   & $\mbi{P} (5,6,8,11) \supset (30)$  & $ YZ^2, W^6, X^5, XY^3 $ & $E_8^2\perp A_1\perp U$ \\
73   & $\mbi{P} (7,8,10,25) \supset (50)$ & $ Z^2, W^6X, X^5Y, Y^5 $ & (19)\\
\hline
\caption{Monomial transformations of the weighted projective spaces.}
\end{longtable}
}

\medskip\noindent
\textsc{Remark.}
We explain the notation of the Table 1.
$1)$ The `No.' follows \cite{Yonemura}. \\
$2)$ We list only the families who have isometric Picard lattices in the famous 95. 
There are many other families of toric $K3$ hypersurfaces. \\
$3)$ In `the vertices of $\Delta$', we state only the vertices of $\Delta$, but not other lattice points in $\Delta$ (e.g. lattice points on edges and faces). For each set of families in Table 1, monomials in the same column (punctuated by commas) correspond. For example, the correspondence between No. 16 and No. 54 shown in Table 1 is determined as follows: 
\[
\begin{array}{ccc}
\textnormal{No.} 16 &   & \textnormal{No.} 54\\
Z^3 & \leftrightarrow & Z^3,\\
W^3YZ & \leftrightarrow & W^3XZ,\\
W^6X & \leftrightarrow & W^7,\\
X^4 & \leftrightarrow & WY^3,\\
WY^3 & \leftrightarrow & X^3Y.
\end{array}
\]
\\
$4)$ `Picard lattices' is due to \cite{Belcastro} and Picard numbers are in the parenthesis. 

\bigskip\noindent
\textsc{Proof.} 
For each $a$, we can choose a polyhedron $\Delta$ as the convex hull of corresponding points in $M(a)$ designated in the table. 
Each birational transform of the weighted projective spaces is given by a correspondence between 
sets of rational monomials. 
It is routine to check that, for weights $a$ and $b$ in each set of rows, 
the correspondence of rational monomials gives an isomorphism between $M(a)$ and $M(b)$, 
and that the polyhedrons are reflexive and isomorphic as lattice polyhedrons. 

$\mbi{P}(a)$, $\mbi{P}(b)$ and $\mbi{P}_\Delta$ contain $\mbi{T}$ in common. 
Thus, for a polynomial $F$ whose Newton polyhedron is in $\Delta$, 
the same zero locus $Z_F$ is contained in those three spaces. 
The compactifications are naturally birational. 
We remark that the zero locus of $F$ in the projective toric varieties may contain some boundary divisors, thus may be different according to the ambient spaces. 
This does not happen for Gorenstein $K3$ hypersurfaces since the boundary divisors are a finite union of toric orbits, thus a finite sum of rational varieties. 
Since $\Delta$ is reflexive, a general anticanonical divisor of $\mbi{P}_\Delta$ has a $K3$ surface as its minimal model. 
Thus there are isomorphisms among the three families of $K3$ surfaces which are minimal models of anticanonical divisors, 
over an open subspace in a projective space, whose dimension is the number of the lattice points in $\Delta$ minus one. 

One can compute the Picard lattice $\Lambda_\Delta$ of a general $K3$ surface by using \cite{Kobayashi}, 
and check that $\Lambda_\Delta$ is isomorphic to both of $\Lambda(a)$ and $\Lambda(b)$. 
For that, it is enough to check that the rank of $\Lambda_\Delta$ coincides with that of $\Lambda(a)$ and $\Lambda(b)$, 
since $\Lambda_\Delta$ contains $\Lambda(a)$ and $\Lambda(b)$, 
and Picard lattices are primitively embedded in the $K3$ lattice by Hodge theory. 
\QED

\vspace{3mm}
For a tetrahedron $\Delta (a)$ with $a$ a weight, let $\textnormal{N}(\Delta (a))$ denote the full Newton polyhedron of $\Delta (a). $

\medskip\noindent
\textsc{Remark.}
If more than three families have the isometric Picard lattice, some subtleties occur. 
For instance, as stated in Table 1, correspondence between No.\,26 and No.\,34 does not fully extend to correspondence including No.\,76. 
One should take a smaller subfamily to establish a correspondence including all three as follows: 

\begin{figure}[htbp]
\begin{center}
\includegraphics[height=.67\linewidth]{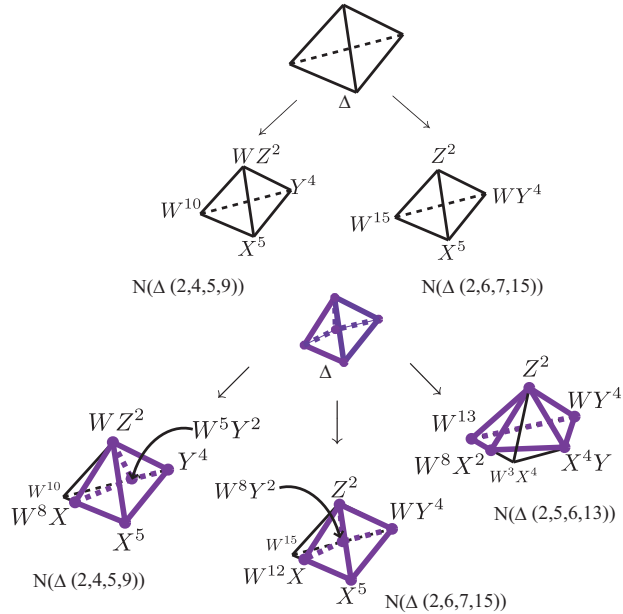}
\end{center}
\caption{above: $\Delta$ for Nos. 26 and 34, \, below: $\Delta$ for Nos. 26, 34 and 76. }
\end{figure}
\noindent
The full Newton polyhedrons of Nos. 26 and 34 are isomorphic, so that $\Delta$ for this pair is to be isomorphic to these polyhedrons, whilst one may have to ``remove" vertices of the full Newton polytopes of Nos. 26, 34 and 76 to obtain $\Delta$ for the set of these three weights. 

\medskip\noindent
\textsc{Remark.}
When $\Delta$ is symmetric, 
clearly other monomial transformations exist; in the list below, the monomials in bold can be exchanged in a row. 
\begin{longtable}{cllc}
\hline
No. & Families & The vertices of $\Delta$ & Picard lattice \\
\hline
16 & $\mbi{P} (3,6,7,8) \supset (24)$ & $ \mbi{Z^3}, W^3YZ, W^6X, X^4, \mbi{WY^3}$ & $E_8 \perp (A_2)^3 \perp U$ \\
54 & $\mbi{P} (3,5,6,7) \supset (21)$ & $\mbi{Z^3}, W^3XZ, W^7, WY^3, \mbi{X^3Y}$ & (16)\\
\hline
30   & $\mbi{P} (5,7,8,20) \supset (40)$ & $ Z^2, W^4Z, \mbi{WX^5}, W^5XY, \mbi{Y^5} $   & $E_8\perp T_{2,5,5}$ \\
86   & $\mbi{P} (4,5,7,9) \supset (25)$  & $ YZ^2, W^4Z, \mbi{X^5}, W^5X, \mbi{WY^3} $  & (18)\\
\hline
46   & $\mbi{P} (5,6,22,33) \supset (66)$ &  $ Z^2, \mbi{W^{12}X}, \mbi{X^{11}}, Y^3 $ & $E_8^2\perp U$ \\
65   & $\mbi{P} (3,5,11,14) \supset (33)$ &  $ XZ^2, \mbi{W^{11}}, \mbi{WX^6}, Y^3 $ & (18)\\
80   & $\mbi{P} (4,5,13,22) \supset (44)$ &  $ Z^2, \mbi{W^{11}}, \mbi{WX^8}, XY^3 $ & \\
\hline
56   & $\mbi{P} (5,6,8,11) \supset (30)$  & $ YZ^2, \mbi{W^6}, \mbi{X^5}, \mbi{XY^3} $ & $E_8^2\perp A_1\perp U$ \\
73   & $\mbi{P} (7,8,10,25) \supset (50)$ & $ Z^2, \mbi{W^6X}, \mbi{X^5Y}, \mbi{Y^5} $ & (19)\\
\hline
\caption{Other monomial transformations.}
\end{longtable}
For example, Nos. 16 and 54 have correspondences as follows:
\begin{figure}[htbp]
\begin{center}
\includegraphics[width=.5\linewidth]{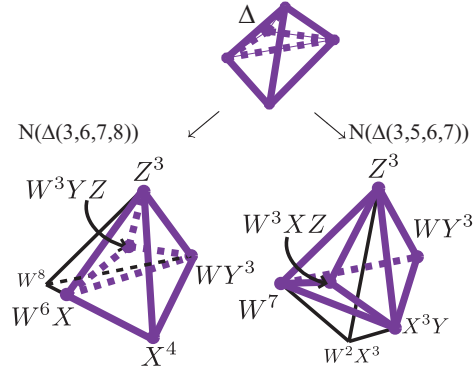}
\end{center}
\caption{Subfamily of Nos. 16 and 54}
\end{figure}

\medskip\noindent
\textsc{Remark.}
The restriction of the Picard group of a resolution of the ambient space $\mbi{P}_\Delta$
 do not always generate $\Lambda_\Delta$. 
We denote by $L_0$ the orthogonal complement of the image of the restriction in the Picard lattice. 

In each set of weights with the isometric Picard lattices, one of the weights has a dual weight system\cite{Ebeling}, 
with one exception, the pair Nos.\,16 and \,54.
There is no reflexive subpolyhedron with $L_0 = 0$ for that pair. 
Although they have a dual weight system, Nos.\,26, 34 and 76, and Nos.\,27 and 49 never contain a reflexive subpolyhedron with $L_0 = 0$. 

\medskip\noindent
\textsc{Remark.}
The real part of logarithmic function gives a homomorphism $(\mbi{C}^{\times})^n \to \mbi{R}^n$; $(z_1,\ldots,z_n) \mapsto (\log |z_1|, \ldots, \log |z_n|)$. 
For a hypersurface $Z$ in $(\mbi{C}^{\times})^n$, 
the image is called the {\it amoeba\/} of $Z$. 

Generally, a monic rational monomial birational map of toric varieties is the morphism which is induced by an isomorphism of the complex tori as complex Lie groups. 
Therefore, it gives a linear isomorphism between amoebas of $K3$ surfaces. 


\end{document}